\documentclass[10pt,reqno]{amsart}
\usepackage[english]{babel}
\usepackage{amsmath}
\numberwithin{equation}{section}
\usepackage{amsfonts}
\usepackage{amssymb}
\usepackage{hyperref}
\hypersetup{citecolor=blue,
    colorlinks=true,
    linkcolor=red,
    filecolor=magenta,      
    urlcolor=cyan,
}
\usepackage{mathrsfs}
\usepackage{amsthm}
\usepackage{graphicx}
\usepackage[usenames]{color}
\usepackage[margin=1in]{geometry}

\renewcommand\hat{\widehat}
\renewcommand\bar{\overline}
\newcommand{\norm}[1]{\left\|#1\right\|}
\newcommand{\abs}[1]{\left|#1\right|}

\renewcommand*{\div}{\ensuremath{\mathrm{div\,}}}

\newcommand*{\R}{\ensuremath{\mathbb{R}}}

\newcommand{\ueq}[1][]{%
  \if\relax\detokenize{#1}\relax
    \sbox0{$\underbrace{=}_{}$}%
    \mathrel{\mathmakebox[\wd0]{=}}
  \else
    \mathrel{\underbrace{=}_{\mathclap{#1}}}
  \fi}

\theoremstyle{plain}
\newtheorem{theorem}{Theorem}[section]

\theoremstyle{definition}
\newtheorem{remark}[theorem]{Remark}

\renewcommand*{\tilde}{\widetilde}
\renewcommand*{\hat}{\widehat}
\newcommand*{\curl}{\ensuremath{\mathrm{curl\,}}}

\def \p{\partial}
\def \intz{\int_{-H_0+b(X)}^{\xi(X,t)}}
\def \intzz{\int_{-1+\beta b }^{\epsilon\xi}}
\def \inti{\int_{z}^{\epsilon\xi}}

\newcommand{\su}{{\vert_{\rm s}}}
\newcommand{\bo}{{\vert_{\rm b}}}
\newcommand{\surf}{{\vert_{z=0}}}
\newcommand{\bott}{{\vert_{z=-1+\beta b}}}
\newcommand{\suadi}{{\vert_{z=\epsilon \xi}}}
\newcommand{\s}{{\vert_{z= \xi}}}

\def \uh {U_{h}}
\def \uv {U_{v}}
\def \uhs {U_{sh}^{\star}}
\def \bh {B_{h}}
\def \bv {B_{v}} 
\def \bhs {B_{sh}^{\star}}
\def \jh {j_{h}}
\def \jv {j_{v}}

\def \U {\mb{U}}

\def \B {\mb{B}}
\def \ome {\mb{\omega}}
\def \j {\mb{j}}

\def \nab {\nabla_{X,z}}

\newcommand{\mb}[1]{\mbox{\boldmath{$#1$}}}

\numberwithin{equation}{section}

\title[Shallow MHD models ]{Asymptotic shallow models arising in magnetohydrodynamics}

\author[D.~Alonso-Or\'an]{Diego Alonso-Or\'an}
\address{University of Bonn, Institute For Applied Mathematics. Endenicher Allee 60, 53115 Bonn, Germany.}
\email{alonso@iam.uni-bonn.de}

\begin{document}
\begin{abstract}
In this paper, we derive a new shallow asymptotic model for the free boundary plasma-vacuum problem governed by the magnetohydrodynamic  equation, vital in describing large-scale processes in flows of astrophysical plasma. More precisely, we present the magnetic analogue of the 2D Green-Naghdi equations for water waves in the presence of weakly shared vorticity and magnetic currents. The method is inspired by developed ideas for hydrodynamics flows in \cite{CastroLannes2} to reduce the $(d+1)$-dimensional dynamics of the problem to a finite cascade of equations which can be closed at the precision of the model. 
\end{abstract}
\maketitle

\section{Introduction}
Plasma is an ionized gas consisting of freely moving positively charged ions, electrons and neutrals. It is by far the most common phase of ordinary matter present in the Universe. We are in constant contact with the small amounts that are not, such as the oceans and seas, but electrically charged fluids are everywhere throughout the  galaxy, \cite{WarDel14}. Macroscopic plasma processes are usually described by the so called magnetohydrodynamic equations (MHD), first proposed by the physicist H. Alfv\'en \cite{Alf42}. Here the behaviour of fluid particles is governed by Newton's second law under the effect of an electromagnetic force described by Maxwell's equations. Furthermore, the movements of present charged particles create an electric field that also affects this magnetic field. 

Motivated by the problem of magnetic plasma confinement in laboratory research, the plasma-vacuum
interface problem has attracted the interest of the mathematical community in the last decades \cite{Tra12,Tra16,GW16,Gu19}. However, studying the full dynamics of the equations is too complex, mainly because the moving surface boundary is part of the solution. This led physicists, oceanographers and mathematicians to derive and replace the original equations by approximate asymptotic systems in specific physical regimes. Those systems are  more amenable to numerical simulations and their properties are more transparent. 

 In the case of the water waves equation describing the motion of an inviscid and incompressible fluid delimited by a free-surface , the most prominent example is the non-linear shallow water equations, also known as Saint-Venant equations, cf. \cite{Ovs, AlvarezLannes}. In the shallow water regime, when $\mu:=\frac{H^2}{L^2} \ll 1$ (where $L$ is the typical horizontal scale and $H$ the typical depth), the non-linear shallow water  system is derived from the free-surface Euler equation by averaging and neglecting the $\mathcal{O}(\mu)$ order terms. As a counterpart  of simplifying the model by dropping the $\mathcal{O}(\mu)$ order terms,  we miss completely the dispersive effects, vital for many applications. Keeping them in the equations  and just neglecting the $\mathcal{O}(\mu^2)$ order terms one obtains the so called Green-Naghdi equations  \cite{ GN76} or Serre equations \cite{Serre53}. We refer to \cite{Makarenko,AlvarezLannes} for a rigorous derivation \footnote{In the sense that their solutions remain close to the exact solution of the free-surface Euler equations of all these models.} and \cite{LannesBonneton} for a more recent review.  Moreover, the Green-Naghdi system is one of the most used model to perform numerical simulations of coastal flows \cite{Cienfuegos,MarcheLannes}. In the equations the flow is assumed to be irrotational or almost irrotational, which breaks down in the presence of rip currents or when underlying currents are present. The complication to describe waves in the presence of non-trivial vorticity is due to its $(d+1)$-dimensional nature while for irrotational flows are only $d$-dimensional. 

Besides different approaches to deal with this difficulty, a novel strategy  was developed in \cite{CastroLannes2}, where the additional terms appearing in the momentum equation can be treated without appealing to the resolution of the $(d+1)$- dimensional vorticity equation. The strategy followed in \cite{CastroLannes2}, inspired by the so called turbulence theory, hinges on deriving a cascade of equations which is actually finite at the precision of the model (($\mathcal{O}(\mu^2)$-order) without any artificial closure. The resulting equations are an extension of the classical Green-Naghdi equations, where the non-hydrostatic pressure terms are affected by the interactions between the horizontal and vertical components of the vorticity.  A rigorous justification of the asymptotic expansions and models derived in \cite{CastroLannes2} have been studied in a companion paper by same authors in \cite{CastroLannes1}. In conclusion, the hydrodynamic shallow water models have been well-studied and extensively used to describe fluid motion in oceans, with direct applications in coastal engineering, \cite{Lannesbook,LannesModelling}.  Surprisingly enough, most of the electrically conducting fluids appearing in astrophysical plasmas, such as accretion disks, planetary atmospheres or stars is that they are, in some sense, thin. However, the role that magnetic field plays when the fluid is electrically conducting is far from being well understood.

The introduction of magnetic effects into the shallow water system was first proposed by Gilman \cite{gilman} and is used as a model of the solar tachocline. The tachocline, first coined by Spiegel and Zahn in\cite{SpiegelZahn}, is a very thin layer in the Sun of thickness of about two and five per-cent of the solar radius which bridged the transition region between the convective zone from the radiative zone, \cite{Hughes}. The approximation of the shallow  MHD equation neglects $\mathcal{O}(\mu)$ order terms and hence can be understood as an magnetic analogue of the non-linear shallow water equations. Since their derivation they have been investigated from a theoretical and modelling point of view.  Studies on linear and non-linear waves  have been treated in \cite{Schecter} as well as applications of the shear-flow instabilities in \cite{Mak}. The shallow MHD system has also been shown to be  hyperbolic and enjoys a Hamiltonian structure \cite{Dellar,Sterck,Rossmanith}.  Extensions from one-layer to  multi-layers models of the shallow MHD equations, experimenting with different stratification settings and enriching the model by choosing different states on each layer (different velocities and different magnetic fields) where also documented in \cite{Hunter,Zeitlin}. \\

The main purpose of this article is to derive new asymptotic models, in dimension $d=1,2$,  for the shallow MHD equation (SMHD) in the presence of vorticity and magnetic currents up to a precision of $\mathcal{O}(\mu^{3/2})$. Our strategy is influenced  by the work in \cite{CastroLannes2} for the water waves equations based on ideas reminiscent of the turbulence theory. The new models take into account the dispersive effects missed in the model by Gilman \cite{gilman} due to the non-hydrostatic pressure and more important manage to deal with the strongly coupling between the vorticity and the current equations, which are assumed to be weakly sheared plasmas (see \S \ref{S:32}). This is achieved, by getting rid of the $(d+1)$-dimensional vorticity-current system, and rather look for an equation involving two $d$-dimensional quantities, i.e.,  the shear velocity and magnetic shear
$$U_{sh}= \int_{z}^{\epsilon \xi} \omega_{h}^{\perp}, \quad  B_{sh}= \int_{z}^{\epsilon \xi} j_{h}^{\perp}, $$
which are the  key terms to close the cascade equation. Roughly speaking, those quantities represent the contributions of the horizontal vorticity and magnetic current into the horizontal momentum equation. We should mention that the asymptotic models derived in this paper are not rigorously justified which is out of the scope of the present article and is left for future work.

\subsection*{Plan of the paper}
The paper is structured as follows. In Section \S \ref{sec:2} we present the basic equations of the free-surface MHD problem. Section \S \ref{sec:3} is devoted to the asymptotic analysis of the free-surface MHD equation, where we first recast the MHD problem using the elevation-discharge formulation. Subsection \S\ref{S:32} deals with non-dimensional averaged-formulation and further remarks on the asymptotic regimes studied on this paper. An asymptotic description of the velocity field and magnetic field is performed in \S \ref{S:32}, which depends on the velocity and magnetic shear term. Section \S \ref{sec:4} presents  the different shallow asymptotic models up to different orders of approximation. In subsection \S \ref{S:41}, we turn to the derivation of the 2D magnetic Green-Naghdi equations, where we first derive an evolution equation for the velocity and magnetic shear \S \ref{S:4.11} and later a closure equation for the tensors in \S \ref{S:4.12}. Lastly, in \S \ref{S:42} we also handle the 1D magnetic Green-Naghdi equation, where several terms trivialize.  Finally a conclusion and future work perspectives are given in Section \S \ref{S:5}.

\subsection*{Notation} We will use the following notation throughout the manuscript.
\begin{itemize}
\item We denote by $d=1,2$ the horizontal dimension, by $X\in \mathbb{R}^d$ the horizontal coordinate and by $z$ the vertical variable. \\
\item For every vector field $\mb{F}\in\mathbb{R}^{3}$ we denote by  $F_{h}$ its horizontal component and by $F_{v}$ its vertical component. If the vector field $\mb{F}\in\mathbb{R}^{2}$ we denote by $\mb{F}^{\perp}=(F_1,F_2)^{\perp}=(-F_2,F_1)$. Let $A\in R^{2}\times R^{2}$ be a matrix, then $A^t$ is the transpose of $A$. Let $A,B\in  \mathbb{R}^{2} \times \mathbb{R}^{2}$ two matrices, we denote by $A:B$ the standard matrices product. For any two different vector fields $\mb{F,G}\in\mathbb{R}^{2}$, we define the tensor product $\mb{F}\otimes \mb{G}= FG^t$, as the usual outer product. \\
\item We use  $\nabla$ to denote  the gradient with respect to the horizontal variables and $\nab$ is the full three-dimensional gradient operator. The rotational and divergence operators are given by 
$$\mbox{curl }\mb{F}=\nab\times \mb{F}, \quad \mbox{div} \mb{F} = \nab\cdot \mb{F}.$$ 
\end{itemize}

\section{The basic equations }\label{sec:2}
\begin{figure}[h]
\includegraphics[width=8cm]{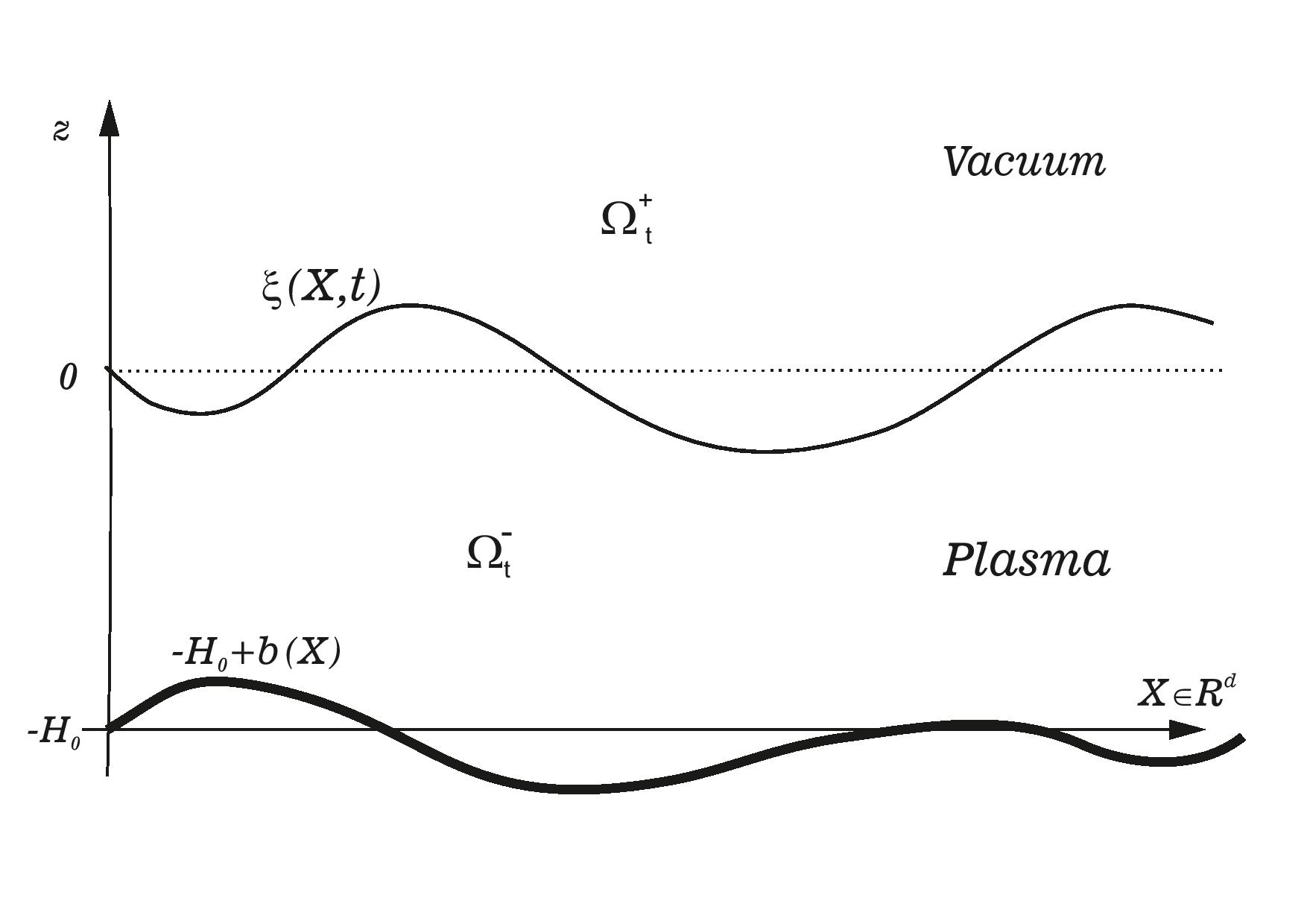}
\centering
\caption{The plasma-vacuum interface problem setting}
\label{fig:1}
\end{figure}

In this section we are concerned with a free-boundary problem of ideal incompressible magnetohydrodynamics. The problem consists in finding a variable domain $\Omega_{t}^{-}$ occupied by an electrically conducting homogeneous plasma, together with a velocity field $\U=\U(X,z,t)$, the scalar pressure $P=P(X,z,t)$ and the magnetic field $\B=\B(X,z,t)$ satisfying the equations of MHD. The elevation of the free surface is parametrized by the graph of a function $\xi(\cdot,t)$, and the non-moving bottom topography is parametrized by a time independent function $-H_{0}+b(X)$ where $H_{0}$ represents the the depth of the plasma and $b$ the possible variation of the bottom, see Figure \ref{fig:1} . Therefore, the domain occupied by the plasma at time $t$ is
  \begin{equation}
\label{plasma:domain}
\Omega_{t}^{-}=\{ (X,z)\in\R^{d+1}: -H_{0} + b(X) < z < \xi(X,t)\}.
\end{equation}
The dynamics of the plasma region is governed by the ideal incompressible MHD equations
\begin{equation}\label{mhd:full}
\left\lbrace
\begin{array}{lll}
\p_{t} \U + (\U\cdot \nab)\U =-\dfrac{1}{\rho}\nab P - g\mb{e_{z}}+\dfrac{1}{\rho\mu_{0}} ((\B\cdot \nab)\B-\dfrac{1}{2}\nab \abs{\B}^{2})\quad \mbox{in} \  \Omega_{t}^{-}, \\  \\
\p_{t} \B+(\U\cdot \nab)\B = (\B\cdot \nab)\,U  \quad  \mbox{in} \  \Omega_{t}^{-}, \\  \\
\div \ \B = 0, \quad \div  \ \U= 0 \quad  \mbox{in} \   \Omega_{t}^{-}.
\end{array}\right.
\end{equation}
where the external forces due to gravity  $\mb{g}=-g \mb{e_{z}}$ are also taken into account. Above $\rho$ is the density assumed constant, and $\mu_{0}$ the magnetic permeability constant. It is assumed that the plasma is surrounded by vacuum region, 
\begin{equation}
\label{vacuum:domain}
\Omega_{t}^{+}=\{ (X,z)\in\R^{2+1}: z > \xi(X,t)\},
\end{equation}
at time $t$. Since, vacuum has no  density, velocity or electric current, the pre-Maxwell dynamics apply \cite{Goed1, Bernstein}. In such a case, the magnetic field $\hat{\B}$ is determined by the \textit{div-curl}
\begin{align}\label{vacuum:eq}
\nab\times \hat{\B}&=0, \quad \mbox{div}\hat{\B}=0, \quad  \mbox{in} \ \Omega_{t}^{+}.
\end{align}
 \\

The plasma–vacuum interface  is now free to move since the plasma is surrounded by vacuum. The physical quantities of the plasma and vacuum region must satisfied non-trivial jump conditions connecting the fields across the interface (cf. \cite{Goed1} for a thorough exposition). The first boundary condition, the so called kinematic boundary condition, is related to the fact that the plasma particles on the surface stay on the surface moving with normal 
$V_{n}=\U\cdot N $ given by
\begin{equation}
\label{kine:boundary}
 \p_{t}\xi-\U_{\su}\cdot N=0
\end{equation}
where $\U(X,t)_{\su}=\U(X,t,\xi(X,t))$ and the $N=(-\nabla\xi, 1)^{T}$. The second and third conditions satisfied at the surface, is the pressure balance condition and magnetic jump continuity
\begin{subequations}
\begin{align}
\left[ \left[ P+\dfrac{1}{2\mu_{0}}\abs{\B}^{2}\right] \right]&=0, \label{press:balance:1} \\
\left[ \left[ \B_{\su}\cdot N \right] \right]&=0,  \label{mag:balance:1}
\end{align}
\end{subequations}
where $\left[\left[ f \right]\right]=\hat{f}-f$ denotes the jump in a quantity across the surface. Hence,
\begin{subequations}
\begin{align}
P+\dfrac{1}{2\mu_{0}}\abs{\B}^{2}&=\dfrac{1}{2\mu_{0}} \abs{\hat{\B}}^{2}, \quad  \text{on} \quad \{ z=\xi(X,t)\},\label{press:balance:2} \\
\B_{\su}\cdot N&=\hat{\B}_{\su}\cdot N \label{mag:balance:2}.
\end{align}
\end{subequations}

Finally, we also impose two boundary conditions at the bottom topography assumed to be a perfect conducting and impermeable material, this is
\begin{subequations}
\begin{align}
\U_{\bo} \cdot N_{b}&=0,\label{bottom:bound:u}\\
 \ \B_{\bo}\cdot N_{b}&=0,\label{bottom:bound:b}
\end{align}
\end{subequations}
where $\U_{\bo}(X,t)=\U(X,t,-H_{0}+b(X))$, $\B_{\bo}(X,t)=\B(X,t,-H_{0}+b(X))$  and $N_{b}=(-\nabla b, 1)^{T}$. \\ \\
In order to simplify the dynamics of the equations we will assume that the vacuum magnetic field $\hat{\B}$ is identically zero, which of course is a trivial solution of \eqref{vacuum:eq}. In this particular case, the free boundary MHD problem is given by
\begin{equation}\label{mhd:full:2}
\left\lbrace
\begin{array}{lll}
\p_{t} \U + (\U\cdot \nab)\U =-\dfrac{1}{\rho}\nab P - g\mb{e_{z}}+\dfrac{1}{\rho\mu_{0}}((\B\cdot \nab)\B-\dfrac{1}{2}\nab \abs{\B}^{2}) \quad \text{in} \  \Omega_{t}^{-}, \\ \\
\p_{t} \B+(\U\cdot \nab)\B = (\B\cdot \nab)\,U  \quad \text{in} \  \Omega_{t}^{-}, \\ \\
\div \ \B = 0, \quad \div  \ \U= 0 \quad \text{in} \  \Omega_{t}^{-}, 
\end{array}\right.
\end{equation}
with boundary conditions
\begin{subequations}
\begin{align}
\p_{t}\xi+\U_{\su}\cdot N&=0, \label{kine:boundary2} \\
(P+\dfrac{1}{2\mu_{0}}\abs{\B}^{2})_{\su}&=0,\label{press:balance:3} \\
\B_{\su}\cdot N&=0, \label{mag:balance:3} \\
\U_{\bo} \cdot N_{b}&=0,\label{bottom:bound:u2}\\
\B_{\bo}\cdot N_{b}&=0.\label{bottom:bound:b2}
\end{align}
\end{subequations}
\begin{remark}
In \eqref{press:balance:3} we are neglecting the effects of the surface tensions, which can be also incorporated into the model, cf. \cite{ChenDing}. Moreover, we also neglect the Coriolis effects induced by planetary rotation. 
\end{remark}
\section{The averaged MHD equations and asymptotic analysis}\label{sec:3}
In this section we will recast the MHD equations \eqref{mhd:full:2} using the so called elevation-discharge formulation (cf. \cite{Lannesbook} )that proves to be very convenient in order to obtain and understand different asymptotic models. The central idea of the formulation is to get rid of the vertical variable  by integrating vertically the  horizontal component of the free surface MHD equation. To that purpose, let us denote by $\U=\U(X,z,t)=(U_{h}(X,z,t),U_{v}(X,z,t))$ , $\B=\B(X,z,t)=(B_{h}(X,z,t),B_{v}(X,z,t))$ the horizontal and vertical component of the plasma velocity and magnetic field, respectively.  Moreover, we introduce the horizontal velocity discharge $Q$,
\begin{equation}\label{discharge:Q}
Q(X,t)=\intz \uh(X,z,t) \  dz,
\end{equation}
and the horizontal magnetic discharge $Q_B$,
\begin{equation}\label{discharge:Qb}
Q_B(X,t)=\intz \bh(X,z,t) \  dz.
\end{equation}
Integrating vertically the horizontal component of equations \eqref{mhd:full:2} and using the boundary conditions \eqref{kine:boundary2}, \eqref{mag:balance:3},\eqref{bottom:bound:u2} and \eqref{bottom:bound:b2} gives
\begin{equation}\label{mhd:averaged1}
\left\lbrace
\begin{array}{lll}
\p_{t} \xi + \nabla\cdot Q=0, \\
\p_{t} Q + \nabla\cdot \left( \displaystyle\int_{-H_{0}+b}^{\xi} U_{h}\otimes U_{h} \right) + \dfrac{1}{\rho} \displaystyle\int_{-H_{0}+b}^{\xi}\nabla P_{B} = \dfrac{1}{\rho\mu_{0}}\nabla\cdot\left(\displaystyle \int_{-H_{0}+b}^{\xi} B_{h}\otimes B_{h} \right), \\
\p_{t} Q_{B} +  \nabla\cdot \left( \displaystyle\int_{-H_{0}+b}^{\xi} B_{h}\otimes U_{h} \right) = \nabla\cdot \left(\displaystyle \int_{-H_{0}+b}^{\xi} U_{h}\otimes B_{h} \right), \\
\nabla\cdot  Q_{B} = 0,
\end{array}\right.
\end{equation}
where $P_B=P+\frac{1}{2\mu_0}\abs{\B}^2$ is the magnetic pressure. Next, we can decompose the pressure magnetic field into a hydrostatic magnetic pressure term $P^h_B$ and non-hydrostatic magnetic pressure term $P_B^{nh}$. It is easy to check that $\xi=0, \U=0, \B=0$ is a particular steady state solution of equations \eqref{mhd:full:2} and therefore the vertical component of the first equation in \eqref{mhd:full:2} gives the following ordinary differential equation for $P_B$:
\begin{equation}\label{pressure:hydrostatic}
-\dfrac{1}{\rho}\p_{z}P_B-g=0
\end{equation}
with boundary condition $P_B\surf=0$, and hence the solution is the hydrostatic magnetic pressure $P^h_{B}=-\rho g z$. Similarly when the fluid is not at rest, integrating the vertical component between $z$ and $\xi$ we have that
$P_B=\rho g(\xi-z) + P^{nh}_{B}$ where the non-hydrostatic magnetic pressure is given by
\begin{equation}\label{pressure:nonhydrostatic}
P^{nh}_{B}(X,z,t)=\rho \int_{z}^{\xi(X,t)} \p_{t}\uh + \U\cdot\nab\uv-\B\cdot\nab\bv \ dz.
\end{equation}
Plugging \eqref{pressure:nonhydrostatic}, we have that the evolution equations are given by
\begin{equation}\label{mhd:averaged2}
\left\lbrace
\begin{array}{lll}
\p_{t} \xi + \nabla\cdot Q=0, \\
\p_{t} Q + \nabla\cdot \left( \displaystyle\int_{-H_{0}+b}^{\xi} U_{h}\otimes U_{h} \right) +gh\nabla\xi+ \dfrac{1}{\rho} \displaystyle\int_{-H_{0}+b}^{\xi}\nabla P^{nh}_{B} = \dfrac{1}{\rho\mu_{0}}\nabla\cdot\left( \displaystyle\int_{-H_{0}+b}^{\xi} B_{h}\otimes B_{h} \right), \\
\p_{t} Q_{B} +  \nabla\cdot \left(\displaystyle \int_{-H_{0}+b}^{\xi} B_{h}\otimes U_{h} \right) = \nabla\cdot \left( \displaystyle\int_{-H_{0}+b}^{\xi} U_{h}\otimes B_{h} \right), \\
\nabla\cdot  Q_{B} = 0,
\end{array}\right.
\end{equation}
where $h(X,t)$ is the plasma height, $h=\xi+H_0+b$. In addition, we also decompose the horizontal velocity and magnetic vector field as 
\begin{subequations}
\begin{align}
\uh(X,t,z)&= \bar{\uh}(X,t)+\uh^{\star}(X,t,z) \label{decomp:u} \\
\bh(X,t,z)&=\bar{\bh}(X,t)+\bh^{\star}(X,t,z) \label{decomp:b}
\end{align}
\end{subequations}
 where for a general function $g(\cdot, t)$ as
\begin{equation}
 \bar{g}(X,t)=\dfrac{1}{h}\intz g(X,t,z) dz 
 \end{equation}
where
\begin{equation}
\quad h(X,t)=H_{0}-b(X)+\xi(X,t) \quad \text{and} \quad g^{\star}(X,t,z)=g(X,t,z)-\bar{g}(X,t). 
\end{equation}
Using \eqref{decomp:u}-\eqref{decomp:b} into the equations, we have that
\begin{equation}\label{mhd:averaged3}
\left\lbrace
\begin{array}{lll}
 \p_{t}\xi+\nabla\cdot Q=0,  \\
\p_{t} Q+ \nabla\cdot (\dfrac{1}{h}Q\otimes Q) + \nabla\cdot\mathcal{R}+gh\nabla \xi+\dfrac{1}{\rho}\displaystyle\intz \nabla P^{nh}_{B} dz =\dfrac{1}{\rho\mu_0}\nabla\cdot (\dfrac{1}{h}Q_B\otimes Q_B)+\dfrac{1}{\rho\mu_{0}}\nabla\cdot \mathcal{R}_{b} , \\
\p_{t}Q_B+ \nabla\cdot (\dfrac{1}{h}Q_b\otimes Q) + \nabla\cdot\mathcal{R}_{m}=\nabla\cdot (\dfrac{1}{h}Q\otimes Q_B) + \nabla\cdot \mathcal{R}^{t}_{m}, \\
\nabla \cdot Q_B=0,
\end{array}\right.
\end{equation}
where 
\[ \mathcal{R}:=\intz \uh^{\star}\otimes\uh^{\star} dz  \quad \text{and} \quad
\mathcal{R}_{b}:=\intz \bh^{\star}\otimes\bh^{\star} dz , \]
\[\mathcal{R}_{m}:=\intz \bh^{\star}\otimes\uh^{\star} dz  \quad \text{and}  \quad  \mathcal{R}^{t}_{m}:= \intz \uh^{\star}\otimes\bh^{\star} dz .\]
We will refer to equations  \eqref{mhd:averaged3} as the  averaged MHD equations. Although the equations are exact, they are not closed. Indeed, quick inspection of equations \eqref{mhd:averaged3} reveal that the terms $\mathcal{R},\mathcal{R}_b,\mathcal{R}_{m}$ and the non-hydrostatic magnetic pressure term are non explicit in terms of the functions $\xi, Q, Q_B$. Actually, as we will see in \S \ref{S:41}, they will depend on the  vorticity $\ome=\nab\times \U$ and magnetic current $\j=\frac{1}{\mu_{0}}\nab\times \B$ equations  given by
\begin{equation}\label{eq:vor:cur}
\left\lbrace
\begin{array}{lll}
\p_{t} \ome + (\U\cdot\nab)\ome-(\B\cdot\nab)\j= (\ome\cdot\nab)\U- (\j\cdot\nab)\B ,  \\ \\
\p_{t} \j + (\U\cdot\nab)\j- (\B\cdot\nab)\ome= (\j\cdot\nab)\U-(\ome\cdot\nab)\B -2 \displaystyle\sum_{i=1}^{3}\nab\U_{i}\times \nab\B_{i}.
\end{array}\right.
\end{equation}
\subsection{Non-dimensional free boundary MHD equation}\label{S:3.1}
To study the behaviour of the solutions of the full system, is generally too complicated, since it contains all the information of the dynamics. We will therefore simplify the terms which seem less important to us. To determine which terms are more relevant, we will use very well known method based on the principle of dimensioning. We non-dimensionality the equations by using several characteristic length of the problem, namely: the typical amplitude of the waves $a_{s}$, the typical depth $H_{0}$, the typical horizontal scale $L$ and the order of bottom variations $a_{b}$. With these characteristic scales we can construct three different dimensionless parameters:
\begin{equation}
\mu = \dfrac{H_{0}^2}{L^2}, \epsilon =\dfrac{a_{s}}{H_{0}}, \beta =  \dfrac{a_{b}}{H_{0}};
\end{equation}
where $\epsilon$ is called the amplitude parameter, $\beta$ the topography parameter and $\mu$ the shallowness parameter. 
\begin{remark}
For simplicity we have assumed here that the horizontal scale $L$ is the same in both the longitudinal $L_x$ and transversal $L_y$ direction, however this could also be considered into the modelling yielding a new parameter $\gamma=\frac{L_x}{L_y}$ called transversality parameter \cite{Lannesbook}. Since the main goal of this article is to derive shallow models we will assume that $\mu\ll 1$, but not other smallness assumption  will be done.
\end{remark}
With the before-mentioned parameters, we will define the following dimensionless variables
\begin{equation}
\label{dimless:variables}
\tilde{z} =\dfrac{z}{H_{0}}, \quad \tilde{X}=\dfrac{X}{L}, \quad \tilde{t}=\dfrac{t}{\dfrac{L}{\sqrt{gH_{0}}}},
\end{equation}
and dimensionless unknowns functions 
\begin{equation}\label{dim:unknown}
\left\lbrace
\begin{array}{lll}
\tilde{\xi}=\dfrac{\xi}{a}, \quad \tilde{\uh}=\dfrac{\uh}{a\sqrt{\dfrac{g}{H_{0}}}}, 	\quad \tilde{\uv}=\dfrac{\uv}{\dfrac{aL}{H_{0}}\sqrt{\dfrac{g}{H_{0}}}}, \\
\tilde{\bh}=\dfrac{\bh}{a\sqrt{\dfrac{\mu\rho g}{H_{0}}}},\quad  \tilde{\bv}=\dfrac{\bv}{\dfrac{aL}{H_{0}}\sqrt{\dfrac{\mu\rho g}{H_{0}}}},  \quad \tilde{P}=\dfrac{P}{\rho g H_{0}}.
\end{array}\right.
\end{equation}

\begin{remark}
The nondimensionalization of the time variable, velocity field and magnetic field are based on a linear analysis, similar to the one performed for the water waves equation (cf. \cite{Lannesbook,CastroLannes1,CastroLannes2}). 
\end{remark}
With this new variables at hand, we denote by
\begin{equation*}
\left\lbrace
\begin{array}{lll}
\U^{\mu}= \begin{pmatrix} \sqrt{\mu}\uh \\ \uv \end{pmatrix}, \quad 
\B^{\mu}= \begin{pmatrix} \sqrt{\mu}\bh \\ \bv \end{pmatrix}, \quad  \\
\nabla^{\mu}= \begin{pmatrix} \sqrt{\mu} \nabla \\ \p_{z} \end{pmatrix} , \quad
\curl^{\mu} = \nabla^{\mu}\times, \quad
 \div^{\mu}=(\nabla^{\mu})^T\cdot,
 \end{array}\right.
\end{equation*}
the free boundary MHD equations \eqref{mhd:full:2} are given in dimensionless form (omitting tildes) by
\begin{equation}\label{mhd:adim}
\left\lbrace
\begin{array}{lll}
\p_{t} \U^{\mu} + \dfrac{\epsilon}{\mu}\U^{\mu}\cdot\nabla^{\mu}\U^{\mu}= \dfrac{1}{\epsilon} \left(\nabla^\mu P_{B^\mu}+1\right) +\dfrac{\epsilon}{\mu}\B^{\mu}\cdot\nabla^{\mu}\B^{\mu}, \\
\p_{t}\B^{\mu} + \dfrac{\epsilon}{\mu}\U^{\mu}\cdot\nabla^{\mu}\B^{\mu}= \dfrac{\epsilon}{\mu}\B^{\mu}\cdot\nabla^{\mu}\U^{\mu}, \\
\div^{\mu}\U^{\mu}=0, \quad \div^{\mu}\B^{\mu}=0,
\end{array}\right.
\end{equation}
posed in the dimensionless fluid domain,
\begin{equation}
\label{dimless:domain}
\Omega(t)^{-}=\{ (X,z)\in\R^{d+1}: -1+\beta b(X)< z < \epsilon \xi(X,t) \}.
\end{equation}
In  the same way, the boundary conditions in dimensionless form are given by
\begin{equation}
\left\lbrace
\begin{array}{lll}
\p_{t}\xi- \dfrac{1}{\mu}\U^{\mu}_\suadi \cdot  N^{\mu}=0,  \\
\B^{\mu}_{\suadi}\cdot N^{\mu}=0,  \\
 {P_{B^{\mu}}}_{\suadi}=0,\\
 \U^{\mu}_{\bott}\cdot N^{\mu}_{b}=0, \\
  \B^{\mu}_{\bott}\cdot N^{\mu}_{b}=0, 
\end{array}\right.
\end{equation}
with $N^{\mu}= \begin{pmatrix} -\epsilon\sqrt{\mu} \nabla \xi \\ 1 \end{pmatrix}$, $N^{\mu}_{b}=\begin{pmatrix} -\beta\sqrt{\mu} \nabla b \\ 1 \end{pmatrix}$.  In a similar way, the dimensionless version of the averaged MHD equations \eqref{mhd:averaged3} are given by
\begin{equation}\label{adim:averaged}
\left\lbrace
\begin{array}{lll}
 \p_{t}\xi+\nabla\cdot Q=0,  \\
\p_{t} Q+ \epsilon\nabla\cdot (\dfrac{1}{h}Q\otimes Q) +\epsilon \nabla\cdot\mathcal{R}+h\nabla \xi+\dfrac{1}{\epsilon}\displaystyle\intzz \nabla P^{nh}_{B} dz =\epsilon\nabla\cdot (\dfrac{1}{h}Q_B\otimes Q_B)+\epsilon \nabla\cdot \mathcal{R}_{b} , \\
\p_{t} Q_B+ \epsilon \nabla\cdot (\dfrac{1}{h}Q_B\otimes Q) +\epsilon \nabla\cdot\mathcal{R}_{m}=\epsilon\nabla\cdot (\dfrac{1}{h}Q\otimes Q_B) + \epsilon\nabla\cdot\mathcal{R}^{t}_{m}, \\
 \nabla \cdot Q_B=0, 
\end{array}\right.
\end{equation}
where the dimensionless tensors and non-hydrostatic term is given by
\begin{eqnarray}
\dfrac{1}{\epsilon}P^{nh}_{B}&=&\int_{z}^{\epsilon \xi}\left( \p_{t} \uv + \epsilon \uh\cdot\nabla \uv+\dfrac{\epsilon}{\mu} \uv\p_{z}\uv-\epsilon \bh\cdot\nabla \bv+\dfrac{\epsilon}{\mu}\bv\p_{z}\bv \right)  \label{adim:npres} \\
\mathcal{R}&=&\intzz \uh^{\star}\otimes\uh^{\star} dz  \quad \text{and} \quad
\mathcal{R}_{b}=\intzz \bh^{\star}\otimes\bh^{\star} dz ,   \label{adim:reynold}  \\
\mathcal{R}_{m}&=&\intzz \bh^{\star}\otimes\uh^{\star} dz \label{adim:mix:reynold} 
\end{eqnarray}

Analogously as in the case the free-surface equation Euler equation in the presence of non-trivial vorticity \cite{CastroLannes1, CastroLannes2}, we require a true comprehension of the behaviour of the vorticity $\ome$  and the current $\j$. To that purpose, let us do the following digression: due to the incompressibility condition,
\begin{equation*}
\left\lbrace
\begin{array}{lll}
\div^{\mu}\U^{\mu}&=&\mu \nabla\cdot \uh+\p_{z}\uv=0, \\
\div^{\mu}\B^{\mu}&=&\mu \nabla\cdot \bh+\p_{z}\bv=0,
\end{array}\right.
\end{equation*}
and the boundary value condition,
\begin{equation*}
\left\lbrace
\begin{array}{lll}
-\beta\mu \nabla b \cdot {\uh}_{\bott} + {\uv}_{\bott}&=0, \\
-\beta\mu \nabla b\cdot {\bh}_{\bott} + {\bv}_{\bott} &=0,
\end{array}\right.
\end{equation*}
we then see that the vertical component $\uv,\bv$ is of order $\mathcal{O}(\mu)$ is the horizontal component   $\uv,\bh$ is of order  $\mathcal{O}(1)$. Moreover, by a straightforward computation we infer that
 \[  \nabla^{\mu}\times\U^{\mu}= \dfrac{H^2}{aL}\sqrt{\dfrac{H}{g}}\ome, \quad  \nabla^{\mu}\times\B^{\mu}= \dfrac{H^2}{aL}\sqrt{\dfrac{H\mu_{0}}{g\rho}}\j,\]
hence we might rescale them as
 \[  \tilde{\ome}= \dfrac{H^2}{aL}\sqrt{\dfrac{H}{g}}\ome, 
\quad\tilde{\j}= \dfrac{H^2}{aL}\sqrt{\dfrac{H\mu_{0}}{g\rho}}\j.\]
The horizontal component, denoted by $\tilde{\omega}_{h}, \tilde{j}_{h}$ respectively, is given by
 \[  \tilde{\omega}_{h}=  \mu \p_{z}\uh^{\perp}- \sqrt{\mu}\nabla^{\perp}\uv, \quad 
\tilde{j}_{h}=  \mu \p_{z}\bh^{\perp}- \sqrt{\mu}\nabla^{\perp}\bv, \]
so $\tilde{\omega}_{h}, \tilde{j}_{h}$ are of order $\mathcal{O}(\mu)$ if we suppose that $\p_{z}\uh^{\perp},\p_{z}\bh^{\perp} $ are of order $\mathcal{O}(1)$. Since we will treat only the case of weakly sheared flows \cite{Teshukov,RG1, RG2,CastroLannes1,CastroLannes2}, this is we will assume that $\tilde{\omega}_{h}, \tilde{j}_{h}$ of order $\mathcal{O}(1)$. To that purpose, we dimensionality the vorticity $\ome$ and current $\j$ as follows
 \[  \tilde{\ome}^{\mu}= \dfrac{L}{a}\sqrt{\dfrac{H}{g}}\ome, 
\quad\tilde{\j}^{\mu}= \dfrac{L}{a}\sqrt{\dfrac{H\mu_{0}}{g\rho}}\j, \]
so that
 \[ \nabla^{\mu}\times\U^{\mu}=\mu \ome^\mu, \nabla^{\mu}\times\B^{\mu}=\mu \j^\mu .\]
where we have dropped the tidal notation. 

\begin{remark}
Deriving different models based on different assumptions on the strength of the vorticity and magnetic current is also possible. However, here we treat the case of weakly shared plasmas, $\ome^\mu, \j^\mu$ are a $\mathcal{O}(1)$ quantity with respect to $\mu$. In \cite{GG, LannesModelling} models with stronger rotational effects  are studied for the water waves problem, however their rigorous well-posedness justification is still an open problem. We recall here that the aim of the present paper is only to derive new asymptotic models for the (SWMHD) equation, but no a rigorous  proof. 
\end{remark}

\subsection{Asymptotic expansion of the inner velocity and magnetic field}\label{S:32}
In this section, we will derive inner asymptotic description of the velocity field $\U^\mu$ and magnetic field $\B^\mu$. First, let us consider the the following boundary value problem satisfied by the velocity field:
\begin{equation}\label{umu:eq}
\left\lbrace
\begin{array}{lll}
\curl^{\mu} \U^{\mu}&=& \mu \ome^{\mu} \quad \mbox{ in }\quad \Omega\\
\div^{\mu} \U^{\mu}&=& 0 \quad \mbox{ in }\quad \Omega \\
\U^{\mu}_{\suadi}\cdot N^{\mu} &=&0 \quad\mbox{ at the bottom}
\end{array}\right.
\end{equation}

The expansion coincides with the one in \cite{CastroLannes2} for the water waves equation (where $\B^\mu=0$) and is given by
\begin{align}
 \uh&:=\bar{\uh}+\uh^{\star}= \bar{\uh}+ \sqrt{\mu}\uhs + \mu \mathbf{L}^{\star}\bar{\uh} + \mu^{3/2} \mathbf{L}^{\star}U_{sh}^{\star}+\mathcal{O}(\mu^2) \label{expansion:uh}, \\
 \uv&= -\mu \left( \nabla \cdot \bar{\uh} (1+z-\beta b)\right)-\mu^{3/2} \nabla\cdot \int_{-1+\beta b}^{z} \uhs + \mathcal{O}(\mu^{2}) \label{expansion:uv} 
 \end{align}
where  the shear velocity $U_{sh}$ is given by
\begin{equation}\label{uhs:shear}
U_{sh}= \int_{z}^{\epsilon \xi} \omega_{h}^{\perp}, \quad \text{where} \quad U^\star_{sh}=\left(\int_{z}^{\epsilon \xi} \omega_{h}^{\perp}\right)^{\star}=\int_{z}^{\epsilon \xi} \omega_{h}^{\perp}-\dfrac{1}{h}\intzz \int_{z}^{\epsilon \xi} \omega_{h}^{\perp}
\end{equation}
and the operator $\mathbf{L}$ acting on a function $f$ is given by
 \begin{equation}\label{operatorL}
  \mathbf{L}f:= \int_{z}^{\epsilon\xi} \nabla\nabla\cdot\int_{-1+\beta b}^{z} f, \quad \mathbf{L}^{\star}f=(\mathbf{L}f)^{\star}.
  \end{equation}
This suggests that the fluctuations of the horizontal velocity $\uh^{\star}=\uh- \bar{\uh}$ is mainly due to the influence of the vorticity contributing at order $\mathcal{O}(\sqrt{\mu})$, whilst the contribution of the vorticity to the vertical component is much smaller appearing at order $\mathcal{O}(\mu^{3/2})$.  \\

Next, we will consider the expansion of the magnetic field  $\B^\mu$, which satisfies:
\begin{equation}\label{bmu:eq}
\left\lbrace
\begin{array}{lll}
\curl^{\mu} \B^{\mu}&=&\mu \j^{\mu} \quad \mbox{ in }\quad \Omega\\
\div^{\mu} \B^{\mu} &=& 0\quad \mbox{ in }\quad \Omega\\
\B^{\mu}_{\suadi}\cdot N^{\mu}&=&0 \quad\mbox{ at the surface}\\
\B^{\mu}_{\suadi}\cdot N^{\mu} &=&0 \quad\mbox{ at the bottom},
\end{array}\right.
\end{equation}
Notice that the boundary value problem for the magnetic field $\B^\mu$, is completely determined by the magnetic current $\j^\mu$. In contrast with the velocity field where the irrotational part of the $\U^\mu$ is determined from the tangential component at the surface, the magnetic field does not have the extra degree of freedom (cf. \textit{Remark 3} in \cite{CastroLannes2}). Taking the horizontal and vertical component of equations we have that
\begin{equation}\label{bmu:comp:eq}
\left\lbrace
\begin{array}{lll}
\sqrt{\mu}\p_{z}\bh-\sqrt{\mu}\nabla \bv-&=& -\mu  \jh^{\perp} \\
\mu \nabla^{\perp}\cdot \bh&=& \mu \jv \\
\mu\nabla\cdot \bh+\p_{z}\bv&=&0\\
-\mu\nabla \xi \cdot \bh+\bv&=&0, \\
-\beta\mu\nabla b \cdot \bh+\bv&=&0
\end{array}\right.
\end{equation}
We plug the following Ansatz 
 $$\bh=\bh^0+\sqrt{\mu}\bh^1+\mu \bh^2+\mu^{3/2}\bh^3+\mathcal{O}(\mu^2),$$
 $$\bv=\bv^0+\sqrt{\mu}\bv^1+\mu \bv^2+\mu^{3/2}\bv^3+\mathcal{O}(\mu^2),$$ 
 and write $\jv=\jv^0+\sqrt{\mu} \jv^1$ to be consistent with the fact that $\div^{\mu} \j^{\mu}= 0$.  We infer that the pair $(\bh^0,\bv^0)$ is determined uniquely be the following system
\begin{equation}
\left\lbrace
\begin{array}{lll}
 \nabla^{\perp}\cdot \bh^0(x,y)&=&\jv^0(x,y), \\
\nabla\cdot(h \bh^0)&=&0,  \\
\bv^0 &\equiv& 0.
\end{array}\right.
\end{equation}

In order to close the asymptotic models up to precision of order $\mathcal{O}(\mu^{3/2})$ (cf. \S \ref{S:42}) we need to do the following mild assumption on the vertical component of the current; assuming that $\jv^0(x,y)=0$. This can be imposed on the initial data, and check that it is conserved during the evolution of the current equation \eqref{eq:vor:cur}, 
\begin{equation}\label{ass:jv0}
\jv^0(x,y)_{\vert_{t= 0}}=0\quad \mbox{implies that} \quad \jv^0(x,y,t)=0 \quad \forall t>0.
\end{equation}
By the second equation in \eqref{bmu:comp:eq}, we can write $$\bh^0(x,y)=\dfrac{1}{h}\nabla^{\perp}\phi^{0}(x,y),$$ for some potential stream function $\phi$ satisfying the the elliptic system
\begin{equation}\label{bo:eq:elliptic}
\left\lbrace
\begin{array}{lll}
\nabla\cdot(\dfrac{1}{h}\nabla\phi^0)=0, \\
\phi^0\to 0, \quad  \mbox{as} \quad  \norm{x}\to\infty
\end{array}\right.
\end{equation}
with $h$ is smooth bounded function. Therefore, invoking standard elliptic estimates we deduce that $\bh^0=0$, \cite{Gilbarg}.\footnote{We are implicitly imposing that we want to deal with physically meaningful functions, i.e.  functions with finite energy $\phi\in L^2(\mathbb{R})$.} In a similar way, we have that for the $\mathcal{O}(\sqrt{\mu})$-order
\begin{equation}\label{b1:eq}
\left\lbrace
\begin{array}{lll}
\bh^1(x,y,z) &=&\bh^1(x,y,\xi)+\displaystyle\int_{z}^{\xi}j^{\perp}_h, \\
\nabla^{\perp}\cdot \bh^1(x,y,z)&=&\jv^1(x,y,z), \\
\nabla\cdot(h \bh^1(x,y,\xi))&=&-  \displaystyle\ \int_{-1}^{\xi}\nabla\cdot\int_{z}^{\xi}j^{\perp}_h,\\
\bv^1 &\equiv &  0.
\end{array}\right.
\end{equation}
Using the third equation we have that
$$\bh^1(x,y,\xi)= -\dfrac{1}{h}\int_{-1}^{\xi}\int_{z}^{\xi}j^{\perp}_h  + \dfrac{1}{h}\nabla^{\perp}\phi^1(x,y),$$ 
for some potential function $\phi^1(x,y)$.
Hence,
\begin{equation}\label{b1:form}
\bh^1(x,y,z)= \int_{z}^{\xi}j^{\perp}_h-\dfrac{1}{h}\int_{-1}^{\xi}\int_{z}^{\xi}j^{\perp}_h + \dfrac{1}{h}\nabla^{\perp}\phi^1(x,y)= \bhs + \dfrac{1}{h}\nabla^{\perp}\phi^1(x,y),
\end{equation}
where the magnetic shear  $\bhs$ is defined as
\begin{equation}\label{bhs:shear}
 \bhs(x,y,z)=  \int_{z}^{\xi}j^{\perp}_h-\dfrac{1}{h}\int_{-1}^{\xi}\int_{z}^{\xi}j^{\perp}_h.
\end{equation}
The unknown potential function $\phi^1(x,y)$ satisfies that
$$ \nabla^{\perp}\cdot(\dfrac{1}{h}\nabla^{\perp}\phi^1)=\jv^1(x,y,z)-\nabla^{\perp}\cdot\bhs . $$
Equations involving the $\mathcal{O}(\mu)$-order terms are given by
\begin{equation}\label{b2eq}
\left\lbrace
\begin{array}{lll}
\nabla^{\perp}\cdot \bh^2(x,y,z)&=& 0,\\
\bh^2(x,y,z)&=&-\bh^2(x,y,\xi)-\displaystyle\int_{z}^{\epsilon\xi}\nabla \bv^2(x,y,z) \ dz ,\\
\bv^2(x,y,z) &= &\beta \nabla b\cdot \bh^0(x,y)-\displaystyle \int_{-1+\beta b}^{z} \nabla\cdot \bh^0(x,y) \ dz, \\
\nabla\cdot(h\bh^2(x,y,\xi))&=&-\displaystyle\int_{-1+\beta b}^{z} \nabla\nabla\cdot\int_{z}^{\epsilon\xi}\bv^2(x,y,z).
\end{array}\right.
\end{equation}
Thus,  since $(\bh^0,\bv^0)=(0,0)$ 
\begin{equation}\label{b2eq:2}
\left\lbrace
\begin{array}{lll}
\nabla^{\perp}\cdot \bh^2(x,y,z)&=& 0,\\
\bh^2(x,y,z)&=& \bh^2(x,y,\xi), \\
\bv^2(x,y,z) &= & 0,\\
\nabla\cdot(h\bh^2(x,y,\xi))&=&0.
\end{array}\right.
\end{equation}
By the last equation in \eqref{b2eq:2} we have that
$$ \bh^2(x,y)=\dfrac{1}{h}\nabla^{\perp}\phi^2(x,y),$$ 
and hence by the first equation in \eqref{b2eq:2}
\begin{equation}\label{b2:eq:elliptic}
\left\lbrace
\begin{array}{lll}
\nabla\cdot(\dfrac{1}{h}\nabla\phi^2)=0, \\
\phi^2\to 0, \quad  \mbox{as} \quad  \norm{x}\to\infty,
\end{array}\right.
\end{equation}
which implies by the argument as in \eqref{bo:eq:elliptic} that $\bh^2=0$.  Finally, the system for the $\mathcal{O}(\mu^{3/2})$-order satisfies
\begin{equation}\label{b3eq}
\left\lbrace
\begin{array}{lll}
\nabla^{\perp}\cdot \bh^3(x,y,z)&=& 0,\\
\bh^3(x,y,z)&=& \bh^3(x,y,\xi)-\displaystyle\int_{z}^{\epsilon\xi}\nabla \bv^3(x,y,z) \ dz, \\
\bv^3(x,y,z) &= & \beta \nabla b\cdot \bh^{1}(x,y,-1+\beta \nabla b)-\displaystyle \int_{-1+\beta b}^{z} \nabla\cdot \bh^1(x,y,z) \ dz, \\
\nabla\cdot(h\bh^3(x,y,\xi))&=&\nabla\cdot \displaystyle\int_{-1+\beta b}^{\epsilon \xi}\int_{z}^{\epsilon\xi} \nabla \bv^3 .
\end{array}\right.
\end{equation}
The vertical component $\bv^3$ is then given by 
\begin{align}\label{b3:vert:form}
\bv^3(x,y,z)&=-\nabla\cdot \displaystyle \int_{-1+\beta b}^{z} \bhs \ ds -\nabla\cdot\left( \dfrac{1}{h} \nabla^{\perp}\phi^1(x,y)(z+1-\beta b)\right) \ dz
\end{align}
Therefore, using \eqref{b1:form} 
\begin{align}\label{b3:form}
\bh^3(x,y,z)&=-\displaystyle\int_{z}^{\epsilon\xi} \nabla\nabla\cdot\int_{-1+\beta b}^{z} \bhs  +\dfrac{1}{h}\displaystyle\int_{-1+\beta b}^{\epsilon \xi}\int_{z}^{\epsilon\xi} \nabla\nabla\cdot\int_{-1+\beta b}^{z} \bhs +\dfrac{1}{h}\nabla^{\perp}\phi^3 \\
&= - \mathbf{L}^{\star}B_{sh}^{\star}+\dfrac{1}{h}\nabla^{\perp}\phi^3(x,y),
\end{align}
where the operator $\mathbf{L}^{\star}$ is defined as in \eqref{operatorL}. By the third equation in \eqref{b3eq} $\phi^3(x,y)$ satisfies
$$ \nabla\cdot(\dfrac{1}{h}\nabla\phi^3)=\nabla^{\perp}\cdot\mathbf{L}^{\star}B_{sh}^{\star}.$$
Therefore, we have that the magnetic field $\B^{\mu}$ can be describe by the following asymptotic expansion
\begin{align}
 \bh &= \sqrt{\mu}\left(\bhs+\dfrac{1}{h}\nabla^{\perp}\phi^1(x,y)\right) + \mu^{3/2}\left( \mathbf{L}^{\star}B_{sh}^{\star}+\nabla^{\perp}\phi^3(x,y)\right)+ \mathcal{O}(\mu^2), \label{expansion:bh} \\
 \bv &= -\mu^{3/2} \nabla\cdot \displaystyle \int_{-1+\beta b}^{z} \bhs \ ds -\mu^{3/2} \nabla\cdot\left( \dfrac{1}{h} \nabla^{\perp}\phi^1(x,y)(z+1-\beta b)\right) + \mathcal{O}(\mu^{2}). \label{expansion:bv} 
\end{align}

\section{The 2D magnetic Green-Naghdi equation }\label{sec:4}
We deal here with the derivation of the magnetic Green-Naghdi type \footnote{We have coined the equation \textit{2D magnetic Green-Naghdi equation} due to the analogy with the Green-Naghdi equations for water waves. However, since we only deal with precision up to $\mathcal{O}(\mu^{3/2})$, this could be understood as a Green-Naghdi equation with medium amplitude wave assumption $\epsilon=\sqrt{\mu}$.} of equation in the two dimensional case ($d=2$), however we will only treat with models up to precision $\mathcal{O}(\mu^{3/2})$.  Structural complications arise in order to close the two dimensioanl cascade of equations if we want to push the expansion further to full order $\mathcal{O}(\mu^{3/2})$. Those difficulties are due to the strong coupling between the magnetic field and velocity field, as we will notice in the computations.

Let us first present the, \textit{nonlinear shallow MHD equations}, which are an approximation of order $\mathcal{O}(\mu)$ of the full MHD equations \eqref{adim:averaged}, in the sense that we drop all the terms of order  $\mathcal{O}(\mu)$. At this order, we do not need to assume that the vertical component $j_v^{0}$ is zero as we did when deriving the asymptotic expansions in \S \ref{S:32}. Actually, the asymptotic expansion would be given by 
\begin{equation}\label{asim:exp:order:1}
\left\lbrace
\begin{array}{lll}
 \bh = \bar{\bh} + \sqrt{\mu}\bhs + \mathcal{O}(\mu),  \\
 \bv =\mathcal{O}(\mu),  \\
 \uh = \bar{\uh}+ \sqrt{\mu}\uhs +\mathcal{O}(\mu) , \\
 \uv= \mathcal{O}(\mu),
\end{array}\right.
\end{equation}
where $\bar{\bh}=\frac{1}{h}(\nabla^{\perp}\phi^{0}(x,y)+\sqrt{\mu}\nabla^{\perp}\phi^{1}(x,y))$. Plugging the asymptotic expansions  into the non-hydrostatic magnetic pressure \eqref{adim:npres} and the tensors \eqref{adim:reynold}-\eqref{adim:mix:reynold} we have that
\begin{equation}\label{pres:rey:mu}
\left\lbrace
\begin{array}{lll}
 P^{nh}_{B}=\displaystyle\int_{z}^{\xi}\left( \p_{t} \uv + \epsilon\uh\cdot\nabla \uv+\dfrac{\epsilon}{\mu} \uv\p_{z}\uv- \epsilon\bh\cdot\nabla \bv+\dfrac{\epsilon}{\mu}\bv\p_{z}\bv \right)=\mathcal{O}(\mu), \\\\
 \mathcal{R}=\mathcal{O}(\mu), \mathcal{R}_{b}=\mathcal{O}(\mu), \mathcal{R}_{m}=\mathcal{O}(\mu)
\end{array}\right.
\end{equation}
The non-hydrostatic and turbulent effects coming from the tensors are not present in this first approximation. Therefore, neglecting terms of order $\mathcal{O}(\mu)$ in \eqref{adim:averaged}, one obtains the following system of equations: 
\begin{equation}\label{nonshallowmhd}
\left\lbrace
\begin{array}{lll}
\p_{t}\xi+\nabla\cdot (h\bar{\uh})=0,  \\ 
\p_{t} (h\bar{\uh})+ \epsilon\nabla\cdot (h\bar{\uh}\otimes\bar{\uh})  = -h\nabla\xi  +  \epsilon\nabla\cdot (h\bar{\bh}\otimes\bar{\bh})\\ 
\p_{t}(h \bar{\bh})+\epsilon\nabla\cdot (h\bar{\uh}\otimes\bar{\bh}) =\epsilon\nabla\cdot (h\bar{\bh}\otimes\bar{\uh}), 
\end{array}\right.
\end{equation}
for $t\geq 0, x\in \mathbb{R}^d$, $h = 1 + \epsilon\xi-\beta b$ and where we have denoted $\bar{\bh}=\frac{1}{h}(\nabla^{\perp}\phi^{0}(x,y)+\sqrt{\mu}\nabla^{\perp}\phi^{1}(x,y))$. The system \eqref{nonshallowmhd} can be understood as the magnetic MHD version of the well-known \textit{nonlinear shallow water equations}, \cite{Ovs}.
\begin{remark}
The \textit{non-linear shallow MHD equations} presented above where first derived by \cite{gilman}, and studied since their derivation intensively, see \cite{Hunter, Mak2} for recent reviews. In \cite{Dellar}, the authors studied the hyperbolic character of the \textit{non-linear shallow MHD equations}. General theory for symmetric hyperbolic systems (cf. \cite{Benzoni}) assures the local well-posedness for initial data $(\xi(x,0),\bar{\uh}(x.0),\bar{\bh}(x.0))\in H^{s}(\mathbb{R}^d)$  with $s>1+\frac{d}{2}$. In recent paper, Trakhinin \cite{Tra19} investigated the structural stability of shock waves and current-vortex sheets in the non-linear shallow MHD equation \eqref{nonshallowmhd}.
\end{remark}

\subsection{The magnetic 2D Green-Naghdi equation}\label{S:41}
Next, we will compute the contributions of the non-hydrostatic terms \eqref{adim:npres} and the tensors \eqref{adim:reynold}-\eqref{adim:mix:reynold} up to a precision of $\mathcal{O}(\mu^{3/2})$. Using the asymptotic expansion \eqref{expansion:uh} and \eqref{expansion:bh} into \eqref{adim:reynold} and \eqref{adim:mix:reynold} yields
\begin{align}
\epsilon \nabla\cdot \mathcal{R}&=\epsilon\nabla\cdot \int_{-1+\beta b}^{\epsilon\xi} \uh^\star \otimes \uh^\star dz =\epsilon \mu  \nabla\cdot \int_{-1+\beta b}^{\epsilon\xi} \uhs\otimes \uhs +\ \mathcal{O}(\epsilon\mu^{3/2}),
\label{operator:R}
\end{align}
\begin{align}
\epsilon \nabla\cdot \mathcal{R}_{b}&=\epsilon\nabla\cdot \int_{-1+\beta b}^{\epsilon\xi} \bh^\star \otimes \bh^\star dz =\epsilon \mu  \nabla\cdot \int_{-1+\beta b}^{\epsilon\xi} \bhs\otimes \bhs +\ \mathcal{O}(\epsilon\mu^{3/2}),
\label{operator:Rb}
\end{align}
\begin{align}
\epsilon \nabla\cdot \mathcal{R}_{m}&=\epsilon\nabla\cdot \int_{-1+\beta b}^{\epsilon\xi} \bh^\star \otimes \uh^\star dz =\epsilon \mu  \nabla\cdot \int_{-1+\beta b}^{\epsilon\xi} \bhs\otimes \uhs +\ \mathcal{O}(\epsilon\mu^{3/2}).
\label{operator:Rm}
\end{align}
Similarly, we can compute the non-hydrostatic magnetic pressure contributions. Dropping $\mathcal{O}(\mu^{3/2})$ terms we get
\begin{align}
\dfrac{1}{\epsilon}\intzz \nabla P^{nh}_{B}&=\dfrac{1}{\epsilon}\intzz \nabla\int_{z}^{\epsilon \xi}\left( \p_{t} \uv + \epsilon \uh\cdot\nabla \uv+\dfrac{\epsilon}{\mu} \uv\p_{z}\uv-\epsilon \bh\cdot\nabla \bv+\dfrac{\epsilon}{\mu}\bv\p_{z}\bv \right) \nonumber \\
&= -\dfrac{\mu}{\epsilon} \intzz \nabla\int_{z}^{\epsilon \xi} \bigg( \p_{t}(\nabla\cdot(\bar{\uh}(1+z-\beta b)))+ \epsilon \bar{\uh}\cdot\nabla(\nabla\cdot(\bar{\uh}(1+z-\beta b))) \nonumber  \\
&\quad \quad +\epsilon(\nabla\cdot(\bar{\uh}(1+z-\beta b)))\p_{z}(\nabla\cdot(\bar{\uh}(1+z-\beta b)))\bigg) + \mathcal{O}(\epsilon\mu^{3/2}) \nonumber \\
&= \mu h \mathfrak{F}(\p_{t}\bar{\uh}+\epsilon\bar{\uh}\nabla\cdot\bar{\uh}) + \mu \epsilon h \mathfrak{D}(\bar{\uh})+ \mathcal{O}(\epsilon\mu^{3/2}) \label{pres:GN}
\end{align}
where the operators $\mathfrak{F}$ and $\mathfrak{D}$ , are defined as 
$$ \mathfrak{F}g=-\dfrac{1}{3h}\nabla(h^3\nabla\cdot g)+\dfrac{1}{2h}(\nabla(h^1 \nabla b\cdot g)- h^2\nabla b g \cdot g)+\nabla b\nabla b\cdot g,$$
$$ \mathfrak{D}(g)= -2\mathfrak{R}_1(\p_{x}g\cdot \p_{y}g^{\perp}+(\nabla\cdot g)^2) + \mathfrak{R}_2(g\cdot (g\cdot \nabla)\nabla b), $$
and 
$$\mathfrak{R}_1 f = -\dfrac{1}{3h}\nabla(h^3f)-\dfrac{h}{2}f\nabla b, \quad \mathfrak{R}_2 f=\dfrac{1}{2h}\nabla (h^2f)+f\nabla b.$$
We have derived the  non-hydrostatic magnetic pressure contributions using the operators $\mathfrak{F}$ and $\mathfrak{D}$ which we borrowed from the formulation in \cite{BCLMT}. Therefore, from inserting the expressions \eqref{operator:R}-\eqref{operator:Rm} and the pressure \eqref{pres:GN} in the equations \eqref{adim:averaged} we infer that 
\begin{equation}\label{2DGN:eq1}
\left\lbrace
\begin{array}{lll}
 \p_{t}\xi+\nabla\cdot (h\bar{\uh})=0,  \\ \\
\p_{t} (h\bar{\uh})+ \epsilon\nabla\cdot (\bar{\uh}\otimes h\bar{\uh})+h\nabla \xi+ +\mu \mathcal{A}_{1} =0, \\ \\
\p_{t} \nabla^{\perp }\phi + \epsilon \nabla\cdot ( \nabla^{\perp}\phi \otimes \bar{\uh}) +\epsilon \sqrt{\mu}\nabla\cdot\mathcal{R}_{m}=\epsilon\nabla\cdot ( \bar{\uh}\otimes \nabla^{\perp}\phi) + \epsilon\sqrt{\mu}\nabla\cdot\mathcal{R}^{t}_{m}, 
\end{array}\right.
\end{equation}
where 
$$ \mathcal{A}_{1}:=  \epsilon \nabla\cdot\mathcal{R}+  h \mathfrak{F}(\p_{t}\bar{\uh}+\epsilon\bar{\uh}\nabla\cdot\bar{\uh})+  \epsilon h \mathfrak{D}(\bar{\uh})-\epsilon \nabla\cdot (\dfrac{1}{h}\nabla^{\perp}\phi\otimes \nabla^{\perp}\phi) - \epsilon \nabla\cdot \mathcal{R}_{b}$$
and
$$ \mathcal{R}=\int_{-1+\beta b}^{\epsilon\xi} \uhs\otimes\uhs \ dz, \quad  \mathcal{R}_{b}=\int_{-1+\beta b}^{\epsilon\xi} \bhs\otimes \bhs \ dz , \quad  \mathcal{R}_{m}=\int_{-1+\beta b}^{\epsilon\xi} \bhs\otimes \uhs \ dz. $$
Noticing that $\p_{t}h + \nabla\cdot(h\bar{\uh})=0$ and denoting by $\bar{\bh}= \sqrt{\mu}\frac{1}{h}\nabla^{\perp }\phi $ we can rewrite equation \eqref{2DGN:eq1} in a more compact way, namely
\begin{equation}\label{2DGN:eq2}
\left\lbrace
\begin{array}{lll}
 \p_{t}\xi+\nabla\cdot (h\bar{\uh})=0,  \\ \\
(1+\mu \mathfrak{F})(\p_{t} \bar{\uh}+ \epsilon \bar{\uh}\cdot\nabla\bar{\uh})+\nabla \xi +\epsilon \mu \mathfrak{D}(\bar{\uh}) + \epsilon \mu \frac{1}{h}\nabla\cdot\mathcal{R}= \epsilon \bar{\bh}\cdot\nabla\bar{\bh} +\epsilon \mu  \frac{1}{h}\nabla\cdot \mathcal{R}_{b} \\ \\
\p_{t} \bar{\bh} + \epsilon \bar{\uh}\cdot\nabla\bar{\bh} +\epsilon \sqrt{\mu}\nabla\cdot\mathcal{R}_{m}=\epsilon \bar{\bh}\cdot\nabla\bar{\uh} + \epsilon\sqrt{\mu}\nabla\cdot\mathcal{R}^{t}_{m}.
\end{array}\right.
\end{equation}
Hence, to close the equations we have find closure equations for the quantities $\mathcal{R},\mathcal{R}_{b},\mathcal{R}_{m}.$  To that purpose we have first to derive an evolution equation $\uhs$  in \eqref{uhs:shear} and magnetic shear $\bhs$ defined in \eqref{bhs:shear}. With this equations at hand, we will be able to find the closure equations and close the system. \\

\subsubsection{\textbf{Evolution equation for $U_{sh}^{\star}$  and $B_{sh}^{\star}$} }\label{S:4.11}
The dimensionless the vorticity and the current are given by
 $$ \ome^\mu = \dfrac{L}{a}\sqrt{\dfrac{H}{g}}\ome, \quad \j^\mu= \dfrac{L}{a}\sqrt{\dfrac{H\mu_{0}}{g\rho}}\j,$$
so that
 \[ \nabla^{\mu}\times\U^{\mu}=\mu \ome^\mu, \quad  \nabla^{\mu}\times\B^{\mu}=\mu \j^\mu .\]
Hence the dimensionless vorticity-current equations of \eqref{eq:vor:cur}
\begin{equation}\label{vor:cur:adi}
\left\lbrace
\begin{array}{lll}
\p_{t} \ome^{\mu} + \dfrac{\epsilon}{\mu}\U^{\mu}\cdot\nabla^{\mu}\ome^{\mu}-\dfrac{\epsilon}{\mu}\B^{\mu}\cdot\nabla^{\mu}\j^{\mu}= \dfrac{\epsilon}{\mu}\ome^{\mu}\cdot\nabla^{\mu}\U^{\mu}- \dfrac{\epsilon}{\mu}\j^{\mu}\cdot\nabla^{\mu}\B^{\mu} ,  \\ \\
\p_{t} \j^{\mu} + \dfrac{\epsilon}{\mu}\U^{\mu}\cdot\nabla^{\mu}\j^{\mu}- \dfrac{\epsilon}{\mu}\B^{\mu}\cdot\nabla^{\mu}\ome^{\mu}= \dfrac{\epsilon}{\mu}\j^{\mu}\cdot\nabla^{\mu}\U^{\mu}-\dfrac{\epsilon}{\mu}\ome^{\mu}\cdot\nabla^{\mu}\B^{\mu} -2\dfrac{\epsilon}{\mu^2} \nabla^{\mu}\U_{i}^{\mu}\times  \nabla^{\mu}\B_{i}^{\mu}.
\end{array}\right.
\end{equation}
Let us first derive an evolution equation for the shear velocity $\uhs$. The computations are an adaptation of \S{2.3} in \cite{CastroLannes2} .We recall the main steps, highlighting the new modifications. The horizontal component of the vorticity equation \eqref{vor:cur:adi}  is given by
\begin{align}
\label{hor:vor:shear}
\p_{t} \omega_{h}^{\mu}+\epsilon \uh\cdot\nabla\omega_{h}^{\mu}+\dfrac{\epsilon}{\mu}\uv\partial_{z}\omega_{h}^{\mu} -\epsilon \bh \cdot \nabla j^{\mu}_{h}-\dfrac{\epsilon}{\mu}\bv\partial_{z}j^{\mu}_{h} \nonumber \nonumber &=\epsilon \omega^{\mu}_{h}\cdot\nabla\uh \nonumber \\
&\quad+ \dfrac{\epsilon}{\sqrt{\mu}}\omega^{\mu}_{v}\partial_{z}\uh  -\epsilon j^{\mu}_{h}\cdot\nabla \bh-\dfrac{\epsilon}{\sqrt{\mu}}j_{v}^{\mu}\partial_{z}\bh
\end{align}
Using the expansion \eqref{expansion:uh}-\eqref{expansion:uv} and \eqref{expansion:bh}-\eqref{expansion:bv}, we have that
\begin{align}
\p_{t}\omega^{\mu}_{h}+\epsilon\bar{\uh}\cdot\nabla\omega_{h}^{\mu}-\epsilon\nabla\cdot[(1+z-\beta b)\bar{\uh}]\p_{z}\omega_{h}^{\mu} &=\epsilon \omega_{h}^{\mu}\cdot\nabla\bar{\uh}  \nonumber \\
& \quad \quad + \epsilon\omega_{v}\p_{z}U^{\star}_{sh}-\epsilon j^{\mu}_{v}\p_{z}B^{\star}_{sh}+ \mathcal{O}(\epsilon\sqrt{\mu}).
\end{align}\label{hor:curr:shear}
Using $\omega^{\mu,\perp}_{h}=-\p_{z}U^{\star}_{sh}$, $j^{\mu,\perp}_{h}=-\p_{z}B^{\star}_{sh}$ and the fact that $\omega^{\mu}_{v}=\nabla ^{\perp}\cdot \bar{\uh}+\mathcal{O}(\sqrt{\mu})$, $j^{\mu}_{v}=j_{v}^{0}+\sqrt{\mu}j_{v}^{1}$  with $j_v^{0}\equiv 0$ as \eqref{ass:jv0} yields
\begin{align}
\p_{t}\omega_{h}^{\mu}&+\epsilon\bar{\uh}\cdot\nabla\omega_{h}^{\mu}-\epsilon\nabla\cdot[(1+z-\beta b)\bar{\uh}]\p_{z}\omega_{h}^{\mu} 
=\epsilon \omega_{h}^{\mu}\cdot\nabla\bar{\uh}- \epsilon(\nabla^{\perp}\cdot\bar{\uh})\omega^{\mu,\perp}_{h}+ \mathcal{O}(\epsilon\sqrt{\mu}).
\end{align}
Taking the perpendicular operator $\perp$, integrating the equation between $z$ and $\epsilon\xi$ and using boundary condition $\p_{t}\xi+\nabla\cdot Q=0$ at the surface
we have that
\begin{align}
\p_{t}U_{sh}+\epsilon\bar{\uh}\cdot\nabla U_{sh}+\epsilon U_{sh}\cdot\nabla\bar{\uh}+\epsilon \left[\nabla\cdot ((1+z-\beta b)\bar{\uh})\omega^{\perp}_{h}\right]=\mathcal{O}(\epsilon\sqrt{\mu}).
\label{eq:ushear}
\end{align}
Applying operator $\displaystyle\frac{1}{h}\int_{-1+\beta b}^{\epsilon \xi}$ to \eqref{eq:ushear} and subtracting the resulting equation from \eqref{eq:ushear}, we infer that
\begin{align}
\label{eq:us:star}
\p_{t}U^{\star}_{sh}+\epsilon\bar{\uh}\cdot\nabla U^{\star}_{sh}+\epsilon U^{\star}_{sh}\cdot\nabla\bar{\uh} - \epsilon \left[\nabla\cdot ((1+z-\beta b)\bar{\uh})\p_{z}U^{\star}_{sh}\right] = \mathcal{O}(\epsilon\sqrt{\mu}).
\end{align}

\begin{remark}
 Due to the assumption $j^{v}_{0}\equiv 0$ in the asymptotic description of the magnetic field described in Subsection \ref{S:32} we notice that equation \eqref{eq:us:star} coincides with the equation \textit{(2.32)} in \cite{CastroLannes2} for water waves. Without the $j^{v}_{0}\equiv 0$ assumption, we could also derive an evolution equation for $\uhs$ involving non-trivial  contributions of the magnetic field, however we would not be able to close the system.
\end{remark}
In a similar way, the horizontal component of the current density equation \eqref{vor:cur:adi}, is given by
\begin{align}
&\p_{t} j^{\mu}_{h}+\epsilon \uh\cdot\nabla j^{\mu}_{h}+\dfrac{\epsilon}{\mu}\uv\partial_{z}j^{\mu}_{h} -\epsilon \bh \cdot \nabla \omega_{h}^{\mu}-\dfrac{\epsilon}{\mu}\bv\partial_{z}\omega_{h}^{\mu} \nonumber \nonumber =\epsilon j^{\mu}_{h}\cdot\nabla\uh + \dfrac{\epsilon}{\sqrt{\mu}}j^{\mu}_{v}\partial_{z}\uh \label{hor:curr:ad:shear} \\
&-\epsilon \omega_{h}^{\mu}\cdot\nabla \bh-\dfrac{\epsilon}{\sqrt{\mu}}j^{\mu}_{v}\partial_{z}\bh-2\frac{\epsilon}{\mu^2}\left(\mu^{3/2}(\nabla^{\perp}\bh\cdot\p_{z}\uh-\nabla^{\perp} \uh\cdot\p_{z}\bh)+\sqrt{\mu}(\nabla^{\perp}\bv\cdot\p_{z}\uv-\nabla^{\perp}\uv\cdot\p_{z}\bv) \right) \nonumber
\end{align}
Plugging in the asymptotic expansion \eqref{expansion:uh}-\eqref{expansion:uv} and \eqref{expansion:bh}-\eqref{expansion:bv} and dropping terms of order $\mathcal{O}(\epsilon\sqrt{\mu})$,
\begin{align}
\p_{t}j^{\mu}_{h}&+\epsilon\bar{\uh}\cdot\nabla j^{\mu}_{h}-\epsilon\nabla\cdot[(1+z-\beta b)\bar{\uh}]\p_{z}j^{\mu}_{h} -\epsilon j^{\mu}_{h}\cdot\nabla\bar{\uh}=- \epsilon(\nabla^{\perp}\cdot\bar{\uh})\p_{z}\bhs+2\nabla^{\perp}\bar{\uh}:\p_{z}\bhs +\mathcal{O}(\epsilon\sqrt{\mu}).
\end{align}
Taking the perpendicular operator $\perp$, integrating in $z$ and using boundary condition $\B^{\mu}_{\s}\cdot N^{\mu}=0$ we have that
\begin{align}
\label{eq:bs}
\p_{t}B_{sh}+\epsilon\bar{\uh}\cdot\nabla B_{sh}-\nabla\cdot\left(\bar{\uh} (1+z-\beta b)\right)\p_{z}B_{sh}+(\nabla\cdot \bar{\uh})B_{sh}  &=- (B^{\perp}_{sh}\cdot \nabla)\bar{\uh}^{\perp}+\nabla^{\perp}\cdot\bar{\uh}B_{sh}^{\perp}  \nonumber\\
&\quad -2\nabla^{\perp}\bar{\uh}^{\perp}:B_{sh}+\mathcal{O}(\epsilon\sqrt{\mu}).
\end{align}
Noticing that $\nabla^{\perp}\bar{\uh}^{\perp}:B_{sh}= (B_{sh}\cdot\nabla^{\perp})\bar{\uh}^{\perp}$ and using the vectorial identity 
$$ (\nabla \cdot F)G+(G^{\perp}\cdot \nabla)F^{\perp}+ (\nabla^{\perp} \cdot F)G^{\perp}=(G\cdot\nabla)A,$$
with $F=\bar{\uh}$ and $G=B_{sh}$ we have that
\begin{align*}
(\nabla \cdot \bar{\uh})B_{sh}+(B_{sh}^{\perp}\cdot \nabla)\bar{\uh}^{\perp}-(\nabla^{\perp} \cdot \bar{\uh})B_{sh}^{\perp}+2(B_{sh}\cdot\nabla) \bar{\uh}^{\perp}=&(B_{sh}\cdot\nabla) \bar{\uh}+2(B_{sh}\cdot\nabla) \bar{\uh}^{\perp} \\
&\quad \quad -2(\nabla^{\perp} \cdot \bar{\uh})B_{sh}^{\perp} \\ 
=&  (B_{sh}\cdot\nabla) \bar{\uh}-2(\nabla^{\perp}\bar{\uh})^{t} B_{sh}^{\perp}.
\end{align*}
Hence, the equation can be rewritten as 
\begin{align}
\label{eq:bs:2}
\p_{t}B_{sh}+\epsilon\bar{\uh}\cdot\nabla B_{sh} +\epsilon(B_{sh}\cdot\nabla) \bar{\uh} -\epsilon\nabla\cdot\left(\bar{\uh} (1+z-\beta b)\right)\p_{z}B_{sh}-2\epsilon(\nabla^{\perp}\bar{\uh})^{t} B_{sh}^{\perp}=\mathcal{O}(\epsilon\sqrt{\mu}).
\end{align}
Taking the average to \eqref{eq:bs:2} and subtracting it from  \eqref{eq:bs:2},
 \begin{align}
\label{eq:bs:3}
\p_{t}\bhs+\epsilon\bar{\uh}\cdot\nabla \bhs +\epsilon(\bhs\cdot\nabla) \bar{\uh} -\epsilon\nabla\cdot\left(\bar{\uh} (1+z-\beta b)\right)\p_{z}\bhs-2\epsilon(\nabla^{\perp}\bar{\uh})^{t} \bhs{^{\perp}}=\mathcal{O}(\sqrt{\mu}).
\end{align}
\subsubsection{\textbf{Closures equation for $\mathcal{R},\mathcal{R}_{b}$ and $\mathcal{R}_{m}$}}\label{S:4.12}
Let us first derive an evolution equation for tensor $\mathcal{R}$. Recall that the evolution equation for the shear velocity $\uhs$ obtained in \eqref{eq:us:star}
$$
\p_{t}U^{\star}_{sh}+\epsilon\bar{\uh}\cdot\nabla U^{\star}_{sh}+\epsilon U^{\star}_{sh}\cdot\nabla\bar{\uh} - \epsilon \left[\nabla\cdot ((1+z-\beta b)\bar{\uh})\p_{z}U^{\star}_{sh}\right] = \mathcal{O}(\epsilon\sqrt{\mu}).
$$
Therefore, taking the time derivative on the tensors $\mathcal{R}$, we have that
\begin{align*}
\p_{t}\mathcal{R} &= \intzz \p_{t}(\uhs\otimes \uhs)+\p_{t}(\epsilon\xi)(\uhs\otimes \uhs)_{\s} \\
&=\displaystyle\sum_{i=1}^{3} J_i+\p_{t}(\epsilon\xi)(\uhs\otimes \uhs)_{\suadi}
\end{align*} 
where
\begin{align*}
 J_1&=  -\epsilon  \intzz \bar{\uh}\cdot\nabla U^{\star}_{sh} \otimes \uhs + \uhs\otimes \bar{\uh}\cdot\nabla U^{\star}_{sh} \\
 &=- \epsilon \bar{\uh}\cdot\nabla\mathcal{R}-\bar{\uh}\cdot\nabla(\beta b)(\uhs\otimes \uhs)_{\bott} + \epsilon\bar{\uh}\cdot\nabla(\epsilon\xi)(\uhs\otimes \uhs)_{\suadi}, 
\end{align*} 

$$
 J_2=  -\epsilon  \intzz  (U^{\star}_{sh}\cdot\nabla)\bar{\uh}  \otimes \uhs + \uhs\otimes ( U^{\star}_{sh}\cdot\nabla)\bar{\uh} =-\nabla\bar{\uh}^{t}\mathcal{R}- \mathcal{R}:\nabla\bar{\uh}
$$
and
\begin{align*}
 J_3&=  \epsilon  \intzz \left[\nabla\cdot ((1+z-\beta b)\bar{\uh})\p_{z}U^{\star}_{sh}\right] \otimes \uhs + \uhs\otimes \left[\nabla\cdot ((1+z-\beta b)\bar{\uh})\p_{z}U^{\star}_{sh}\right] \\
 &= \epsilon  \intzz \nabla\cdot ((1+z-\beta b)\bar{\uh})(\p_{z}U^{\star}_{sh} \otimes \uhs + \uhs\otimes \p_{z}U^{\star}_{sh}) \\
 &= h (\uhs\otimes \uhs)_{\suadi}-\epsilon \mathcal{R}+\nabla(\beta b)\cdot\bar{\uh} (\uhs\otimes \uhs)_{\bott}-\nabla(\beta b)\cdot\bar{\uh} (\uhs\otimes \uhs)_{\suadi}.
\end{align*}
Using the fact that $\p_t \xi + \nabla\cdot (h\bar{\uh})=0$, we obtain that
$$\p_{t} \mathcal{R} +\epsilon\bar{\uh}\cdot\nabla \mathcal{R}+\epsilon(\nabla\cdot\bar{\uh}) \mathcal{R}+\epsilon \nabla\bar{\uh}^{t} \mathcal{R}+\epsilon  \mathcal{R}:\nabla\bar{\uh}=0. $$

Mimicking the computations but for the magnetic shear equation for $\bhs$ given by
$$
\p_{t}\bhs+\epsilon\bar{\uh}\cdot\nabla \bhs +\epsilon(\bhs\cdot\nabla) \bar{\uh} -\epsilon\nabla\cdot\left(\bar{\uh} (1+z-\beta b)\right)\p_{z}\bhs-2\epsilon(\nabla^{\perp}\bar{\uh})^{t} \bhs{^{\perp}}=\mathcal{O}(\sqrt{\mu}),
$$
we arrive to 
$$ \p_{t} \mathcal{R}_{b} +\epsilon\bar{\uh}\cdot\nabla \mathcal{R}_{b}+\epsilon(\nabla\cdot\bar{\uh}) \mathcal{R}_{b}+\epsilon \nabla\bar{\uh}^{t} \mathcal{T}+\epsilon  \mathcal{R}_{b}:\nabla\bar{\uh}-\epsilon(\nabla^{\perp}\bar{\uh})^{t} \mathcal{R}^{S}_{b}=0,$$
where 
$$\mathcal{R}^{S}_{b}=\intzz {\bhs}^{\perp}\otimes \bhs + \bhs\otimes {\bhs}^{\perp} \ dz. $$
Finally we compute the closure equation for $\mathcal{R}_{m}$, which differentiating the tensor $\mathcal{R}_{m}$ gives
\begin{align*}
\p_{t}\mathcal{R}_{m} &= \intzz \p_{t}(\bhs\otimes \uhs)+\p_{t}(\epsilon\xi)(\bhs\otimes \uhs)_{\s} \\
&=\displaystyle\sum_{i=1}^{4} K_i+\p_{t}(\epsilon\xi)(\bhs\otimes \uhs)_{\suadi}
\end{align*} 
with
\begin{align*}
 K_1&=  -\epsilon  \intzz (\bar{\uh}\cdot\nabla) B^{\star}_{sh} \otimes \uhs + \bhs\otimes (\bar{\uh}\cdot\nabla) U^{\star}_{sh} \\
 &=- \epsilon \bar{\uh}\cdot\nabla\mathcal{R}_{b}-\bar{\uh}\cdot\nabla(\beta b)(\bhs\otimes \uhs)_{\bott} + \epsilon\bar{\uh}\cdot\nabla(\epsilon\xi)(\bhs\otimes \uhs)_{\suadi}, 
\end{align*} 
and
$$
 K_2=  -\epsilon  \intzz  (B^{\star}_{sh}\cdot\nabla)\bar{\uh}  \otimes \uhs + \bhs\otimes ( U^{\star}_{sh}\cdot\nabla)\bar{\uh} =- \nabla\bar{\uh}\mathcal{R}_{b}- \mathcal{R}_{b}:\nabla\bar{\uh}^{t},
$$
while
\begin{align*}
 K_3&=  \epsilon  \intzz \left[\nabla\cdot ((1+z-\beta b)\bar{\uh})\p_{z}B^{\star}_{sh}\right] \otimes \uhs + \bhs\otimes \left[\nabla\cdot ((1+z-\beta b)\bar{\uh})\p_{z}U^{\star}_{sh}\right] \\
 &= \epsilon  \intzz \nabla\cdot ((1+z-\beta b)\bar{\uh})(\p_{z}B^{\star}_{sh} \otimes \uhs + \bhs\otimes \p_{z}U^{\star}_{sh}) \\
 &= h (\bhs\otimes \uhs)_{\suadi}-\epsilon\mathcal{R}_{b}+\nabla(\beta b)\cdot\bar{\uh} (\bhs\otimes \uhs)_{\bott}-\nabla(\beta b)\cdot\bar{\uh} (\bhs\otimes \uhs)_{\suadi}.
\end{align*}
To last term is given by
$$
 K_4=  2\epsilon  \intzz (\nabla^{\perp}\bar{\uh})^t  (B^{\star}_{sh})^{\perp}\otimes \uhs =-2\epsilon(\nabla^{\perp}\bar{\uh}^{\perp})^t \mathcal{R}_{m}.
$$
Collecting all the computations, and recalling that $\p_t \xi + \nabla\cdot (h\bar{\uh})=0$, we have that
$$ \p_{t} \mathcal{R}_{m} +\epsilon\bar{\uh}\cdot\nabla \mathcal{R}_{m}+\epsilon(\nabla\cdot\bar{\uh}) \mathcal{R}_{m}+\epsilon \mathcal{R}_{m}\nabla\bar{\uh}^{t}+\epsilon \nabla\bar{\uh}:\mathcal{R}_{m}+2\epsilon(\nabla^{\perp}\bar{\uh}^{\perp})^{t}\mathcal{R}_{m}=0.$$
\subsubsection{Full 2D magnetic Green-Naghdi equation}
We can now express the \textit{two-dimensional magnetic Green-Naghdi equation}, dropping $\mathcal{O}(\mu^{3/2})$ terms by
\begin{equation}\label{2DGN:eq3}
\left\lbrace
\begin{array}{lll}
 \p_{t}\xi+\nabla\cdot (h\bar{\uh})=0,  \\ \\
(1+\mu \mathfrak{F})(\p_{t} \bar{\uh}+ \epsilon \bar{\uh}\cdot\nabla\bar{\uh})+\nabla \xi +\epsilon \mu \mathfrak{D}(\bar{\uh}) + \epsilon \mu \frac{1}{h}\nabla\cdot\mathcal{R}= \epsilon \bar{\bh}\cdot\nabla\bar{\bh} +\epsilon \mu  \frac{1}{h}\nabla\cdot \mathcal{R}_{b} \\ \\
\p_{t} \bar{\bh} + \epsilon \bar{\uh}\cdot\nabla\bar{\bh} +\epsilon \sqrt{\mu}\nabla\cdot\mathcal{R}_{m}=\epsilon \bar{\bh}\cdot\nabla\bar{\uh} + \epsilon\sqrt{\mu}\nabla\cdot\mathcal{R}^{t}_{m}. \\ \\
\p_{t} \mathcal{R} +\epsilon\bar{\uh}\cdot\nabla \mathcal{R}+\epsilon(\nabla\cdot\bar{\uh}) \mathcal{R}+\epsilon \nabla\bar{\uh}^{t} \mathcal{R}+\epsilon  \mathcal{R}:\nabla\bar{\uh}=0. \\ \\
\p_{t} \mathcal{R}_{b} +\epsilon\bar{\uh}\cdot\nabla \mathcal{R}_{b}+\epsilon(\nabla\cdot\bar{\uh}) \mathcal{R}_{b}+\epsilon \nabla\bar{\uh}^{t} \mathcal{R}_{b}+\epsilon  \mathcal{R}_{b}:\nabla\bar{\uh}-\epsilon(\nabla^{\perp}\bar{\uh})^{t} \mathcal{R}^{S}_{b}=0, \\ \\
 \p_{t} \mathcal{R}_{m} +\epsilon\bar{\uh}\cdot\nabla \mathcal{R}_{m}+\epsilon(\nabla\cdot\bar{\uh}) \mathcal{R}_{m}+\epsilon \mathcal{R}_{m}\nabla\bar{\uh}^{t}+\epsilon \nabla\bar{\uh}:\mathcal{R}_{m}+2\epsilon(\nabla^{\perp}\bar{\uh}^{\perp})^{t}\mathcal{R}_{m}=0
 \end{array}\right.
\end{equation}
where the operators $\mathfrak{F}, \mathfrak{D}$ are defined as 
$$ \mathfrak{F}g=-\dfrac{1}{3h}\nabla(h^3\nabla\cdot g)+\dfrac{1}{2h}(\nabla(h^1 \nabla b\cdot g)- h^2\nabla b g \cdot g)+\nabla b\nabla b\cdot g,$$
$$ \mathfrak{D}(g)= -2\mathfrak{R}_1(\p_{x}g\cdot \p_{y}g^{\perp}+(\nabla\cdot g)^2) + \mathfrak{R}_2(g\cdot (g\cdot \nabla)\nabla b), $$
and 
$$\mathfrak{R}_1 f = -\dfrac{1}{3h}\nabla(h^3f)-\dfrac{h}{2}f\nabla b, \quad \mathfrak{R}_2 f=\dfrac{1}{2h}\nabla (h^2f)+f\nabla b,$$
while the tensors $\mathcal{R}, \mathcal{R}_{b}, \mathcal{R}_{m}$ and $\mathcal{R}^{S}_{b}$ stand for
$$ \mathcal{R}=\int_{-1+\beta b}^{\epsilon\xi} \uhs\otimes\uhs \ dz, \quad  \mathcal{R}_{b}=\int_{-1+\beta b}^{\epsilon\xi} \bhs\otimes \bhs \ dz , \quad  \mathcal{R}_{m}=\int_{-1+\beta b}^{\epsilon\xi} \bhs\otimes \uhs \ dz., $$
$$\mathcal{R}^{S}_{b}=\intzz {\bhs}^{\perp}\otimes \bhs + \bhs\otimes {\bhs}^{\perp} \ dz=M \mathcal{R}_{b}+  \mathcal{R}_{b} M^{t}, \quad \mbox{with} \  M=\begin{pmatrix}
  0 & -1\\ 
  1 & 0
\end{pmatrix}.$$

\begin{remark}
If we compare the equations with the one derived for the rotational water-waves in \cite{CastroLannes2}, we observe new two new phenomena. First we have a self interaction of the magnetic shear $B_{sh}$ coded in $\mathcal{R}_{b}$ and the interactions between the shear velocity and magnetic shear gathered in the term $\mathcal{R}_{m}$ and $\mathcal{R}^{t}_{m}$. 
\end{remark}

\subsection{The 1D magnetic Green-Naghdi equation}\label{S:42}
For the sake of completeness we also derive the equations in the one-dimensional case. We consider velocity and magnetic fields given by
\begin{equation*}
\U^{\mu}= \begin{pmatrix} \sqrt{\mu}u \\ 0 \\ \uv \end{pmatrix}, \quad 
\B^{\mu}= \begin{pmatrix} \sqrt{\mu}b \\  0 \\ \bv \end{pmatrix}. \quad  
\end{equation*}
Therefore the scalar vorticity and magnetic current are given by 
\begin{equation*}
\ome^{\mu}= \begin{pmatrix} 0 \\ \omega^{\mu}(x,z,t) \\ 0 \end{pmatrix}, \quad 
\j^{\mu}= \begin{pmatrix} 0 \\  j^{\mu} (x,z,t) \\ 0 \end{pmatrix}. \quad  
\end{equation*}
Moreover we will have that $\uhs=(u^{\star}_{sh},0)^t$ and $\bhs=(b^{\star}_{sh},0)^t$, and hence
$$u^{\star}_{sh}=-\left( \inti \omega^{\mu}(x,z,t) \ dz \right)^{\star}, \quad b^{\star}_{sh}=-\left( \inti j^{\mu}(x,z,t) \ dz\right)^{\star}.$$
Therefore, straightforward modifications of the two-dimensional case, yield the following \textit{one-dimensional magnetic Green-Naghdi equations}

\begin{equation}\label{1DGN:eq}
\left\lbrace
\begin{array}{lll}
 \p_{t}\xi+\p_{x}(h\bar{u})=0,  \\ \\
(1+\mu \mathfrak{F})(\p_{t} \bar{u}+ \epsilon \bar{u}\p_{x}\bar{u})+\p_{x}\xi +\epsilon \mu \mathfrak{D}\bar{u} + \epsilon \mu \frac{1}{h}\p_{x}\mathcal{R}= \epsilon \mu  \frac{1}{h}\p_{x}\mathcal{R}_{b} \\ \\
\p_{t} \mathcal{R} +\epsilon\bar{u}\p_{x} \mathcal{R}+3\epsilon\p_{x}\bar{u} \mathcal{R}=0. \\ \\
\p_{t} \mathcal{R}_{b} +\epsilon\bar{u}\p_{x}\mathcal{R}_{b}+3\epsilon\p_{x}\bar{u}\mathcal{R}_{b}=0, 
 \end{array}\right.
\end{equation}
with the one-dimensional versions of  $\mathfrak{F}, \mathfrak{D}$ are given by
$$ \mathfrak{F}g=-\dfrac{1}{3h}\p_{x}(h^3\p_{x} g)+\dfrac{1}{2h}(\p_{x}(h^2 g \p_{x}b)-h^2\p_{x}g\p_{x}b)+g(\p_{x}b)^2,$$
$$ \mathfrak{D}f=-\dfrac{2}{3h}\p_{x}(h^3(\p_{x}f)^2)+h(\p_{x}f)^2\p_{x}b+\dfrac{1}{2h}\p_{x}(h^2f^2\p^2_{x}b)+f^2\p^2_{x}b\p_{x}b  $$
while we recall the one dimensional tensors $\mathcal{R}$ and $\mathcal{R}_{b}$ are given by 
$$ \mathcal{R}=\int_{-1+\beta b}^{\epsilon\xi} \abs{u^{\star}_{sh} }^2\ dz, \quad  \mathcal{R}_{b}=\int_{-1+\beta b}^{\epsilon\xi} \abs{b^{\star}_{sh} }^2\ dz. $$
\begin{remark}
We notice that the main difference between the one-dimensional setting and the two-dimension version of the \textit{magnetic Green-Naghdi equations} is that the one dimensional version of the averaged MHD equations \eqref{mhd:averaged3} trivialize to 
\begin{equation}\label{1d:mhd:averaged}
\left\lbrace
\begin{array}{lll}
 \p_{t}\xi+\p_{x}Q=0,  \\ 
\p_{t} Q+ \p_{x} (\dfrac{1}{h}\abs{Q}^2) + \p_{x}\mathcal{R}+gh\p_{x} \xi+\dfrac{1}{\rho}\displaystyle\intz \p_{x} P^{nh}_{B} dz =\dfrac{1}{\rho\mu_{0}}\p_{x} \mathcal{R}_{b} , \\ 
\p_{t}Q_B=0 \\ 
\p_{x} Q_B=0.
\end{array}\right.
\end{equation}
Hence the magnetic elevation discharge $Q_{B}=0$, if we want to $Q_{B}$ has finite energy, therefore excluding the interaction of the velocity and magnetic shear (present in \eqref{2DGN:eq3} as $\mathcal{R}_{m}$) and the coupling between the magnetic field $\bar{\bh}$ in the momentum equation, also present in \eqref{2DGN:eq3}.
\end{remark}
\section{Conclusion and future work}\label{S:5}
In this paper, we have derived a new shallow water models for the free-surface MHD equation in the presence of vorticity and magnetic currents. The most essential ingredient of the model is that the resulting equations are $d$-dimensional which avoid the 
need to solve the $(d+1)$-dimensional nature of the vorticity-current equations, reducing the complexity of the system from a mathematical and numerical point of view.  The strategy follows closely the ideas developed in \cite{CastroLannes2}. It is shown that the additional terms appearing in the momentum equation due to the presence of the vorticity and current effects, satisfy additional two dimensional advection type equations which couple to the system. The advected quantities describe the self-interactions of the  shear velocity and magnetic shear induced by the vorticity-current system, and the coupled interactions between the shear velocity and magnetic shear.  \\

Different techniques were developed to perform numerical simulations for the non-linear shallow MHD equation first derived in \cite{gilman}, for example, the constrained transport approach \cite{Sterck2,Rossmanith2}, projection method \cite{Rossmanith} or central upwind scheme \cite{KurTad, Zia}. Therefore a natural perspective is to take into account the vorticity-current effects in the numerical simulations, allowing the modelling of underlying currents in plasmas or sheared plasmas, using the models derived in this article. Another future research direction, is related to the rigorous justification or the well-posedness of the derived models, which to the best of the authors knowledge is an open problem even for the non-linear shallow MHD equations derived by \cite{gilman}.


\vspace{.2in}
\noindent{\bf{Acknowledgment.}} 
The author is deeply indebted to \'{A}ngel Castro  for many useful discussions and suggestions which have significantly improved this manuscript. He also acknowledges helpful conversations with Daniel Faraco and David Lannes. The author has been supported by the ICMAT Severo Ochoa project SEV-2015-0554 grant, MTM2017-85934-C3-2-P, ERC grant 834728 Quamap and the Alexander von Humboldt Foundation.

\newcommand{\etalchar}[1]{$^{#1}$}

\end{document}